\documentclass{article}
\usepackage{graphicx}
\usepackage{amssymb,amsmath}
\usepackage[spanish]{babel}

\def\timehm{\count31=\time \count32=\count31 \divide\count31 by 60
\number\count31 \multiply\count31 by 60 \advance\count32 by
-\count31 :\ifnum\count32<10 0\fi \number\count32}

\newtheorem{pro}{Proposici\'on}[section]

\def\Q{\mathbb{Q}}

\newcommand{\Frob}{{\rm Frob }}

\newcommand{\PSL}{{\rm PSL}}

\newcommand{\PGL}{{\rm PGL}}

\newcommand{\F}{{\mathbb{F}}}

\begin{document}

\title{{\bf A control theorem for the images of Galois actions on certain infinite families of modular forms
}}
\author{Luis Dieulefait
\\
Dept. d'\'{A}lgebra i Geometria, Universitat de Barcelona\\
e-mail: ldieulefait@ub.edu\\
 }
\date{\empty}

\maketitle

\vskip -20mm
%

%

\section{A letter with the result}

$\qquad \qquad \qquad \qquad \qquad \qquad \quad \quad \qquad \qquad$       November 11, 2006\\

Dear Colleagues,\\

   After reading a recent preprint by G. Wiese on the images of Galois representations attached to classical modular forms and applications to inverse Galois theory I have been thinking again on the results of Ribet and Serre giving ``image as large as possible for almost every prime" in the non-CM case. The result of Wiese is used to realize as Galois group over $\Q$ the group $\PGL_2(\ell^r)$ or $\PSL_2(\ell^r)$ for a fixed $\ell$ and exponent $r$ larger than any given exponent $r_0$. I haven't read his proof in detail, but since it uses ``good-dihedral" primes  and the result is for fixed $\ell$, I can imagine that both the good-dihedral prime and the modular form are constructed ad hoc to realize only the desired linear group in this specific characteristic (and large-but-unknown exponent). In any case, since his method works for every $\ell$, the result that he obtains is quite interesting.\\

On the other hand, of course, you and I believe that a stronger result is true, namely, a uniform result: so, assume for a moment that I just consider, for a fixed small prime $t$ and every exponent $n$ a modular form $f_n$ of level $u_n \cdot t^n$ and weight $2$ (and trivial character) (*), where $u_n$ is prime to $t$, without CM. Assume also that we know somehow that the images of the $\ell$-adic and mod $\ell$ Galois representations attached to $f_n$ are ``as large as possible" for every prime $\ell$ and for every $f_n$. In other words, the special family $f_n$ has image large FOR EVERY PRIME instead of just for almost every prime (which is what Ribet's result shows in general for any modular form without CM). 
If such a family $f_n$ exists, the well-known relation between the conductor and the minimal field of definition for a Galois representation with values on a finite field (as explained in Serre's Duke 1987 paper on modularity conjectures and exploited by Brumer to show that for a  large power $t^n$ in the level the field of coefficients of the projective representations attached to $f_n$ must contain the real part of a cyclotomic field of $t^m$-roots of unity, with $m$ going to infinity as $n$ does) implies that the family $f_n$ gives another proof of Wiese's Galois realizations result, and in a uniform way: for any given $\ell$ and exponent $r_0$ we know that taking any element $f_n$ with $n$ sufficiently large in our family we will be realizing not only the desired projective linear group in characteristic $\ell$ and with exponent greater than $r_0$ but also a similar linear group for a set of primes of density as close to $1$ as desired (but smaller than $1$),  always with exponent larger than $r_0$. So, instead of realizing the desired ``linear group over a finite field with large exponent" for an isolated prime our modular forms $f_n$ will do the job for a large density set of primes. \\

As far as I know, it is not yet known how to construct an infinite family of modular forms with growing level $f_n$ as the one described above, having all of them large image for every prime. But I can construct a family with a slightly  weaker property, that is still good enough to derive the above conclusion regarding realizations of linear groups as Galois groups over $\Q$. In particular the family $f_n$ that I have found with levels as in (*) has, of course, the property that the degree of the corresponding field of coefficients goes to infinity with $n$ (because of Brumer result and the factor $t^n$ in the level), and concerning the images of the corresponding Galois representations it has the property that for each $f_n$ with $n \geq4$ we can give an upper bound to the set of exceptional primes computed as a function only of the level of $f_n$ (i.e., all the information we need is the value of the level, not a single eigenvalue is needed) and in particular we can easily show that for any given prime $\ell>3$ there is a value $n_0$ such that $\ell$ is not exceptional for $f_n$ for any exponent $n>n_0$, where here ``exceptional" means dihedral, reducible or (for $\ell = 5$) some of the other cases of small image in Dickson's result.\\

Let us show one example of such a family $f_n$: for any $n \geq 4$ take $f_n$ to be ANY modular form of level $2\cdot 3^n$, weight $2$, trivial character. Because of semistability at $2$ none of them has CM. Since the large ramification at $3$ for $n$ sufficiently large makes easy to see  (again, using the ideas of Serre and Brumer on conductors) that the small special groups in Dickson's list can not occur  for $\ell = 5$ and it is well known that these groups can not occur for larger $\ell$ in weight $2$, we can concentrate in the two problems that have to be solved to control the images of the representations attached to $f_n$ for any $n$: to control dihedral primes and to control reducible primes. In both cases we will use the large ramification at $3$ and the semistability at $2$ to do so.\\

Dihedral primes: Let $n\geq 4$ and $\ell>3$ be a dihedral prime for $f_n$. Using the arguments created by Serre and Ribet, we see that the only possibility is that the mod $\ell$ representation is induced from the quadratic number field ramifying only at $3$. Also, since a dihedral image does not contain unipotent elements, this mod $\ell$ representation is unramified at $2$. Since $2$ is a non-square mod $3$ this implies that the trace $a_2$ of the image of $\Frob \; 2$ in this mod $\ell$ representation has to be $0$ (i.e.: as usual in the residually CM case, half of the traces have to be $0$, and $a_2$ is in that half). On the other hand the $\ell$-adic representation attached to $f$ has semistable ramification at $2$, so we are in a case of raising the level (or lowering the level, depending on the perspective), and as you know very well this can only happen if $a_2 = \pm 3$ (these are numbers which only exist mod $\ell$). Putting the two things together we conclude that $0$ and $\pm 3$ are the same mod $\ell$, and since $\ell>3$ this gives us a contradiction.\\

Reducible primes: this time assume $\ell>3$ is a reducible prime for $f_n$, $n\geq 4$. We anticipate that now such a prime can exist (for example $7$ is reducible for some newforms of level $162$) but we just want to bound the set of reducible primes in terms of the level $2\cdot3^n$. Again, we will use the local information at $2$ and $3$ to do so.
For simplicity of the exposition, we assume that $n= 2\cdot u$ is even. Since the mod $\ell$ representation is reducible (we semisimplify if necessary, so assume it is semisimple) and using the value of the level of $f_n$ (and, because $\ell>3$, it is well-known that residually the conductor at $3$ will be exactly $3^n$) we know that it is just the direct sum $\chi\cdot\psi \oplus \psi^{-1}$, where $\chi$ is the mod $\ell$ cyclotomic character and $\psi$ is a character  of conductor exactly $3^u$ (remember that $u$ is half of $n$, so it is at least $2$). Computing the trace of the image of $\Frob \; 2$ for the mod $\ell$ representation this time we obtain $a_2 = 2\cdot \psi(2) +   \psi^{-1}(2)$. Here the important thing to observe is that the order of $\psi$ is $\phi(3^u) = 2\cdot3^{u-1}$ (where $\phi$ is Euler's function), and that $2$ is primitive modulo $3^u$, so the order of the element $\psi(2)$ is also $2\cdot3^{u-1}$. On the other hand, using again raising the level since the $\ell$-adic representation is semistable at $2$   we must have $a_2 = \pm 3$. Comparing the two formulas for $a_2$ the first observation is that the roots of the characteristic polynomial of the image of $\Frob \; 2$ are $1, 2$ or $-1,-2$, in particular they belong to the prime field $\F_\ell$, so $\psi(2)$ must be in this field, and looking at the order of this element this means that $3^{u-1}$ divides $\ell-1$ (@). This already shows that for $n$ sufficiently large any prime $\ell$ given a priori will not be reducible (thus, will not be exceptional), because the maximal power of $3$ dividing $\ell-1$ is finite. \\
Just for fun, let us bound the set of possibly reducible primes for $f_n$: Comparing the two formulas for $a_2$ (comparing the roots of the polynomials deduced from both formulas) and using the information on the order of $\psi(2)$ we conclude that any reducible $\ell$ must satisfy, in addition to (@), the condition: $\ell$ divides $2^{2\cdot3^{(u-1)}} - 1$. So, this is a bound for the set of reducible primes for $f_n$. For example for $n=4$ (thus $u=2$) we conclude that $\ell$ has to be congruent to $1$ mod $3$ and divides $2^6-1 = 63$, thus the prime $7$ may be reducible (and it is so for some newforms of level $2\cdot 3^4 = 162$), but it is the only possible reducible prime $\ell > 3$ (computing reducible primes for all newforms of this level using the method in my thesis confirms this fact).\\

Conclusion: For any newform in the family $f_n$ described above,  if $n = 2 \cdot u$ (we assume it is even for simplicity)   the residual image is ``large" for any prime which is not congruent to $1$ modulo $3^{u-1}$. This, together with the fact proved by Serre and Brumer that as $n$ (the $3$-part of the conductor) goes to infinity the exponents of the fields of coefficients of the  projective residual representations also go to infinity, has as a corollary that with our family $f_n$ we are realizing, for any prime $\ell>3$, projective linear groups over the field of $\ell^r$ elements  for $r$ arbitrarily large as Galois groups over $\Q$, and we are realizing these groups in a uniform way (i.e., for sufficiently large $n$ we obtain these groups not only for a given $\ell$ but also for large density sets of primes, all with large exponent). Each of them is realized as an extension unramified outside $6\cdot \ell$.\\

Of course we can construct other similar examples taking other suitable pairs of primes instead of $2$ and $3$, we can also take more general levels having semistable ramification at more primes, and other variations. The main point is that we can bound the set of exceptional primes for ALL modular forms in an infinite family of increasing conductor, which is an interesting result that of course can not be obtained using just the computational method explained in my thesis years ago. It is interesting to observe how a very simple ramification condition (one semistable prime and other dividing the conductor with a large power)  was enough to obtain ``uniformly large" images, to explain that only primes that are ``very split" can be exceptional, thus to generate a lot of large exponent linear groups as Galois groups. Maybe other combinations of ramification conditions can lead to similar, or even stronger, results.\\

The idea used to control dihedral primes is an idea I had in Paris in 2002, when considering dihedral primes for the case of Q-curves coming from diophantine equations. The new idea is the idea to control reducible primes, which I had in Berkeley last week during the modularity conference (but I knew since I saw Wiese's paper months ago that the arguments of Serre and Brumer should be key: obtaining  linear groups over fields with large exponents as Galois groups using these results is something I wanted to do already when starting my thesis).  \\

Best regards,\\

               $\qquad \qquad \qquad $               Luis Dieulefait

\end{document}